\documentclass[leqno,12pt]{amsart}
\usepackage{amssymb}

\calclayout

\newtheorem{theorem}{Theorem}

\newtheorem{lemma}{Lemma}

\newcommand{\cF}{{\mathcal F}}
\newcommand{\cG}{{\mathcal G}}
\newcommand{\cI}{{\mathcal I}}
\newcommand{\cL}{{\mathcal L}}
\newcommand{\cO}{{\mathcal O}}
\newcommand{\cQ}{{\mathcal Q}}
\newcommand{\cV}{{\mathcal V}}
\newcommand{\cY}{{\mathcal Y}}

\newcommand{\mG}{{\mathbb G}}
\newcommand{\mP}{{\mathbb P}}

\newcommand{\tZ}{{\tilde Z}}

\newcommand{\pt}{\text{\it pt}}
\newcommand{\Fl}{\text{\it Fl}}
\newcommand{\Hom}{\text{\it Hom}}
\newcommand{\Ext}{\text{\it Ext}}

\newcommand{\codim}{\operatorname{codim}}
\newcommand{\Chow}{\operatorname{Chow}}
\newcommand{\SO}{\operatorname{SO}}
\newcommand{\Supp}{\operatorname{Supp}}

\title{Multiplicity--free subvarieties of flag varieties}

\author{Michel~Brion}
\address{Universit\'e de Grenoble I\\
D\'epartement de Math\'ematiques\\
Institut Fourier, UMR 5582 du CNRS\\
38402 Saint-Martin d'H\`eres Cedex, France}
\email{Michel.Brion@ujf-grenoble.fr}

\begin{document}

\begin{abstract}
Consider a flag variety $\Fl$ over an algebraically closed field, 
and a subvariety $V$ whose cycle class is a multiplicity--free sum
of Schubert cycles. We show that $V$ is arithmetically normal and
Cohen--Macaulay, in the projective embedding given by any ample
invertible sheaf on $\Fl$. 
\end{abstract}

\maketitle

\section*{Introduction}

In this note, we obtain geometric properties of a class of
subvarieties of flag varieties, that includes Schubert varieties.
By a flag variety, we mean a projective homogeneous space $\Fl=G/P$
where $G$ is a semisimple algebraic group and $P$ is a parabolic
subgroup. Then $\Fl$ admits a cellular decomposition into Bruhat
cells with closures the Schubert varieties. Any subvariety 
$V\subseteq Fl$ is rationally equivalent to a linear combination of
Schubert cycles with uniquely determined integer coefficients; these
are in fact non--negative. We say that $V$ is {\it multiplicity--free}
if these coefficients are $0$ or $1$; in loose words, the numerical
invariants of $V$ are as small as possible. 

\smallskip

Clearly, Schubert varieties are multiplicity--free.
Further examples include: the diagonal in $\Fl\times\Fl$ (and, more
generally, those subvarieties of $\Fl\times\Fl$~that are invariant
under the diagonal $G$--action), the irreducible hyperplane sections
of $\Fl$ in its smallest projective embedding, and also the
irreducible hyperplane sections of Schubert varieties in
Grassmannians, embedded by their Pl\"ucker embedding.

\smallskip

From these examples, we see that multiplicity--free subvarieties may
be singular. But their singularities are well--behaved, as shown by
our main result.

\begin{theorem}\label{cone}
Any multiplicity--free subvariety $V$ of a flag variety $\Fl$ is
normal and Cohen--Macaulay. Moreover, $V$ admits a flat degeneration
in $\Fl$ to a reduced Cohen--Macaulay union of Schubert varieties.

As a consequence, for any globally generated invertible sheaf $\cL$ on
$\Fl$, the restriction map $H^0(\Fl,\cL)\to H^0(V,\cL)$ is surjective,
and $H^n(V,\cL)=0$ for any $n\ge 1$. If, in addition, $\cL$ is ample,
then 
$H^n(V,\cL^{-1})=0$ for any $n\le \dim(V)-1$.

Thus, $V$ is arithmetically normal and Cohen--Macaulay in the
projective embedding given by any ample invertible sheaf on $\Fl$.
\end{theorem}

The latter result is well--known in the case of Schubert varieties,
and several proofs are available, based on methods of Frobenius
splitting \cite{Rn85} or standard monomial theory \cite{Ln98}; the
normality of Schubert varieties has also been obtained by Seshadri
\cite{Si87}, using an argument of descending induction. This line of
argument was taken up in \cite{Bn01} to obtain a version of Theorem
\ref{cone} concerning orbit closures of spherical subgroups of $G$,
in characteristic zero ([loc.cit.], Theorem 6 and Corollary 8). The 
latter assumption plays an essential role there, since the proof
relies on the existence of equivariant desingularizations and on the
notion of rational singularities. 

\smallskip

In the present note, we show how to adapt the arguments of \cite{Bn01}
\S 4 to all multiplicity--free subvarieties, and to arbitrary
characteristics. In Section 1, we formulate a criterion 
for normality and Cohen--Macaulayness of certain subvarieties of
arbitrary $G$--varieties (Theorem 2). Section 2 presents some preliminary
results on linearized sheaves that are needed for the proof of our
criterion, given in Section 3. The final Section 4 presents a proof of
Theorem 1, based on that criterion combined with known vanishing
theorems for line bundles on unions of Schubert varieties.

\smallskip

This raises the question of characterizing the cycle classes
of multiplicity--free subvarieties among all multiplicity--free sums
of Schubert cycles. More generally, characterizing the cycle
classes of (irreducible) subvarieties of flag varieties is a very
natural problem. 

\smallskip

For classes of curves, a complete solution follows from the (much more
precise) results of \cite{Pn02}; in particular, each positive class
of one--cycles arises from an irreducible curve. But this 
does not extend to higher dimensions: consider, for example, an
irreducible surface $V\subset \mP^n\times\mP^n$. Decompose the cycle
class of $V$ as  
$x[\mP^2\times\pt] + y [\mP^1\times\mP^1] + z [\pt\times\mP^2]$
where $x,y,z$ are non--negative integers. 
Then $y^2 - xz \geq 0$, as follows from the Hodge index theorem. In
particular, the multiplicity--free class 
$[\mP^2\times\pt] + [\pt\times\mP^2]$ does not
arise from an irreducible surface (alternatively, one may observe
that $(\mP^2\times\pt)\cup(\pt\times\mP^2)$ is not Cohen--Macaulay,
and apply Theorem \ref{cone}.) But any other multiplicity--free class
arises from an irreducible surface, since 
$[\mP^2\times\pt] + [\mP^1\times\mP^1]$ is the class of the blow--up
of a point in $\mP^2$, and 
$[\mP^2\times\pt] + [\mP^1\times\mP^1] + [\pt\times\mP^2]$
is the class of the diagonal in $\mP^2\times\mP^2$.

\smallskip

The problem of characterizing the cycle classes of {\it nonsingular}
subvarieties of flag varieties is also quite natural and open; here
very interesting recent results are those of \cite{Bt00}, for
Grassmannians over the field of complex numbers.

\section{A criterion for normality and Cohen--Macaulayness of
subvarieties of $G$--varieties}

We begin by fixing notation and recalling some results on linear
algebraic groups; a general reference is \cite{Sr98}.

The ground field $k$ is algebraically closed, of arbitrary
characteristic. A {\it variety} is a separated integral scheme of 
finite type over $k$; by a {\it subvariety}, we mean a closed
subvariety. Given a linear algebraic group $\Gamma$, a 
{\it $\Gamma$--variety} is a variety $Z$ endowed with a 
$\Gamma$--action
$$
\sigma:\Gamma\times Z\to Z,~(g,x)\mapsto g\cdot x.
$$

We will mostly consider $P$--varieties where $P$ is a
parabolic subgroup of a semisimple algebraic group $G$, that is, 
$P$ contains a Borel subgroup $B$ of $G$. Then $P$ is connected, 
and the homogeneous space $G/P$ is projective; the quotient map 
$G\to G/P$ is a locally trivial principal $P$--bundle. The number of
$B$--orbits in the flag variety $G/P$ is finite, and each of them 
(a {\it Bruhat cell}) is isomorphic to an affine space; the orbit
closures are called {\it Schubert varieties}.

We fix a Borel subgroup $B\subseteq G$ and consider only 
parabolic subgroups $P\supseteq B$; these are called 
{\it standard}. The quotient map $P\to P/B$ is a locally trivial 
principal $B$--bundle, and $P/B$ is a flag variety for a maximal
semisimple subgroup of $P$.

A standard parabolic subgroup $P$ is {\it minimal}, if $P\ne B$ and
$P$ is minimal for this property; equivalently, $P/B$ is isomorphic to
the projective line $\mP^1$. Then $B$ acts on the cotangent space to
$P/B$ at the base point $B/B$ by a character $\alpha$, the 
{\it simple root} associated with $P$; this sets up a bijective
correspondence between minimal parabolic subgroups and simple
roots. The group $G$ is generated by its minimal parabolic subgroups.

Let $Y$ be a $B$--variety and $P$ a standard parabolic subgroup; we
review the construction of the {\it induced variety} $P\times^B Y$. 
The group $B$ acts on the product $P\times Y$ by 
$b\cdot (p,y)=(pb^{-1},by)$. Since $P\to P/B$ is a locally trivial
principal $B$--bundle, then so is the quotient
$$
q:P\times Y \rightarrow (P\times Y)/B=P\times^B Y,
$$
and $P\times^B Y$ is a variety. The action of $P$ on itself by left
multiplication makes $P\times^B Y$ into a $P$--variety; the projection
$$
p_1:P\times Y\rightarrow P
$$ 
yields a $P$--equivariant map
$$
f:P\times^B Y\rightarrow P/B,
$$ 
a locally trivial morphism. Its fiber at the base point $B/B$ is
$B\times^B Y\cong Y$. This defines a closed $B$--equivariant immersion 
$$
i:Y\rightarrow P\times^B Y.
$$

Next we assume that $Y$ is a $B$--invariant subvariety of a 
$G$--variety $X$. Then the map 
$P\times Y\to X,~(p,y)\mapsto p\cdot y$ 
factors through a $P$--equivariant morphism
$$
\pi:P\times^B Y\rightarrow X
$$
with image $P\cdot Y$. One checks that the product map
$$
f\times\pi:P\times^B Y\rightarrow P/B\times X,~(p,y)B\mapsto (pB,py)
$$ 
is a closed immersion. Thus, $\pi$ is proper (so that $P\cdot Y$ is
a $P$--subvariety of $X$; in particular, $G\cdot Y$ is a
$G$--subvariety), and its fibers identify with closed subschemes of
$P/B$. 

The {\it normalizer} of $Y$ in $G$ is 
$N_G(Y)=\{g\in G~\vert~ g\cdot Y =  Y\}$, a standard parabolic 
subgroup of $G$. For any standard parabolic subgroup 
$P\subseteq N_G(Y)$, the morphism $f\times\pi$ yields an isomorphism 
$$
P\times^B Y\cong P/B\times Y.
$$

We say that a minimal parabolic subgroup $P$ {\it raises} $Y$ if
it is not contained in $N_G(Y)$, that is, $Y$ is strictly contained
in $P\cdot Y$. Then we have $\dim(P\cdot Y)=\dim(Y)+1$, and the
morphism $P\times^B Y\to P\cdot Y$ is generically finite.

If $Y$ is not $G$--invariant, then $Y$ is raised by some minimal
parabolic subgroup. As a consequence, there exists a sequence
$(P_1,\ldots,P_c)$ of such subgroups, so that 
$G\cdot Y=P_1\cdots P_c \cdot Y$, where $c$ denotes 
the codimension of $Y$ in $G\cdot Y$. We say that the $B$--subvariety
$Y$ is {\it multiplicity--free}, if the morphisms 
$$
P_i\times^B P_{i+1}\cdots P_c\cdot Y\rightarrow 
P_i \cdot P_{i+1} \cdots P_c\cdot Y
$$ 
are birational for all such sequences $(P_1,\ldots,P_c)$. Equivalently,
the compositions
$$
P_1\times^B P_2\times^B\cdots \times^B P_c\times^B Y \rightarrow 
G\cdot Y,~ (p_1,p_2,\ldots,p_c,y)B^c\mapsto p_1 p_2\cdots p_c\cdot y
$$
are birational. 

If $G$ is simply--laced (i.e., its simple factors have
type $A$, $D$ or $E$) and $k$ has characteristic zero, then it
suffices to check this condition for one such sequence, by \cite{Bn01}
Proposition 5. But this does not extend to arbitrary $G$, as shown by
the following geometric variant of \cite{Bn01} Example 4.

\medskip

\noindent
{\bf Example.} Let $q$ be a nondegenerate quadratic form on the vector
space $k^5$ (where the characteristic of $k$ is not equal to $2$). The
corresponding special orthogonal group $G=\SO(5)$ acts on the
projective space $\mP^4=\mP(k^5)$ and leaves invariant the nonsingular
quadric $\cQ=\mP(q=0)$, of dimension $3$. Choose a point $p_0\in\cQ$
and a line $\ell_0$ in $\mP^4$, contained in $\cQ$ and passing through 
$p_0$. Let $B$ be the isotropy group in $G$ of the pair
$(p_0,\ell_0)$; this is a Borel subgroup of $G$. The nontrivial
standard parabolic subgroups are the isotropy groups $P_1$, $P_2$ of
$\ell_0$ resp.~$p_0$; both are minimal.

Consider the diagonal action of $G$ on $X=\cQ\times\cQ$, and the
subset $Y=\ell_0\times (\cQ\cap p_0^{\perp})$, where $p_0^{\perp}$
denotes the hyperplane in $\mP^4$, orthogonal to $p_0$ with respect to
$q$. Clearly, $Y$ is a $B$--invariant subvariety of codimension $3$ in
$X$. One checks that 
$P_1 \cdot P_2 \cdot P_1 \cdot Y = X = P_2 \cdot P_1 \cdot P_2\cdot Y$
and that the morphism 
$P_1 \times^B P_2 \times^B P_1 \times^B Y \rightarrow X$
is birational, whereas the morphism
$P_2 \times^B P_1 \times^B P_2 \times^B Y \rightarrow X$
has degree $2$.

\medskip

Returning to the general situation, note that the $B$--subvariety $Y$
is multiplicity--free if and only if

\begin{enumerate}

\item
The morphism $\pi_Y:P\times^B Y\to P\cdot Y$ is birational for any
minimal parabolic subgroup $P$ raising $Y$, and 

\item
the $B$--subvariety $P\cdot Y$ is multiplicity--free.

\end{enumerate}

Note also that any Schubert variety $S$ in a flag variety $G/Q$ is
multiplicity--free. Indeed, $S$ decomposes as a
product $P_1\cdots P_d \cdot Q/Q$, where $(P_1,\ldots,P_d)$ is a
sequence of minimal parabolic subgroups, and  $d=\dim(S)$. Moreover,
$\pi_S$ is birational for any minimal parabolic subgroup raising
$S$. As a consequence, $S$ is multiplicity--free, and the map 
$P_1\times^B \cdots \times^B P_d\times^B Q/Q\to S$ is birational, with 
source a nonsingular projective variety. This morphism is called a 
{\it standard resolution} of the Schubert variety $S$.

The notion of a multiplicity--free $B$--subvariety generalizes
that of a multiplicity--free subvariety of a flag variety, as shown by 

\begin{lemma}\label{multfree}
Let $\Fl =G/Q$ be a flag variety, $V\subseteq \Fl$ a subvariety, and
put 
$$
Y=\{g\in G~\vert~ g^{-1}Q/Q\in V\}.
$$ 
Then $Y$ is a subvariety of $G$, invariant under the $Q$--action by
left multiplication. Moreover, $V$ is normal (resp.~Cohen--Macaulay,
multiplicity--free) if and only if the $B$--subvariety $Y$ is.
\end{lemma}

\begin{proof}
By definition, $Y$ is the preimage of $V$ under the map 
$G\to G/Q,~g\mapsto g^{-1}Q$. Since that map is a locally trivial
principal bundle for the connected group $Q$, then $Y$ is a
subvariety, and it is normal (resp.~Cohen--Macaulay) if and only if
$V$ is. 

To show the equivalence between both notions of
multiplicity--freeness, we note that $V$ is multiplicity--free if and
only if its pull--back to $G/B$ under the natural map $G/B\to G/Q$
is (since the pull-back of any Schubert variety in $G/Q$ is a Schubert
variety in $G/B$). Moreover, $V$ and its pull--back to $G/B$ define 
the same subvariety $Y\subseteq G$. Thus, we may assume that $Q=B$. 

Let now $P$ be a minimal parabolic subgroup and consider the natural
map 
$$
f:G/B\to G/P,
$$ a locally trivial fibration with fiber $\mP^1$.
Then the subvariety of $G$ associated with $f(V)\subseteq G/P$ (or,
equivalently, with $f^{-1}f(V)\subseteq G/B$) is $P\cdot Y$; as a
consequence, $P$ raises $Y$ if and only if the restriction 
$$
f_V:V\to f(V)
$$ 
is generically finite. Moreover, the fibers of $f_V$ identify to those
of $\pi_Y:P\times^B Y\to P\cdot Y$; in particular, both morphisms have
the same degree $d_V$. Note that $d_V\ne 0$ if and only if $\pi_Y$ (or,
equivalently, $f_V$) is birational. 

Let $S\subseteq G/B$ be a Schubert variety of positive dimension. 
Write $S=P_1\cdots P_d/B$, where $P_1,\ldots,P_d$ are minimal parabolic
subgroups and $d=\dim(S)$; let $P=P_d$ and $S'=P_1\cdots P_{d-1}/B$.
Then $S=f^{-1}f(S)$, and the restriction $S'\to f(S')=f(S)$ is
birational. We thus obtain the equalities of intersection numbers: 
$$\displaylines{
\int_{G/B} [V]\cdot [S] = \int_{G/B} [V]\cdot f^*[f(S)]
=\int_{G/P} f_*[V] \cdot [f(S)] 
\hfill\cr\hfill
= d_V \int_{G/P} [f(V)]\cdot f_*[S']
=d_V \int_{G/B} [f^{-1}f(V)] \cdot [S'],
\cr}$$
as follows from the projection formula and from the equalities
$f_*[V]=d(V)[f(V)]$, $f_*[S']=[f(S)]$. Moreover, any Schubert variety
of positive codimension in $G/B$ arises as some $S'$. Thus, $V$ is
multiplicity--free if and only if: $d_V=1$ and $f^{-1}f(V)$ is
multiplicity--free, for any minimal parabolic subgroup $P$ raising
$Y$. This completes the proof.
\end{proof}

We may now state our criterion for normality and
Cohen--Macaulayness. It may be regarded as a characteristic--free
version of \cite{Bn01} Theorem 5 (see the remarks at the end of
Section 3 for more detailed comments).

\begin{theorem}\label{main}
Let $X$ be a $G$--variety and $Y\subseteq X$ a multiplicity--free
$B$--subvariety containing no $G$--orbit.

If $G\cdot Y$ is normal, then $Y$ is normal as well.

If, in addition, $G\cdot Y$ is Cohen--Macaulay, then $Y$ is
Cohen--Macaulay as well.
\end{theorem}

\section{Induction of linearized sheaves}

In this section, we study an analogue of the induction of varieties,
concerning linearized sheaves; these may be defined as follows, after
\cite{Kf78}.

Let $\Gamma$ be a linear algebraic group, $Z$ a $\Gamma$--variety with
action map $\sigma:\Gamma\times Z\rightarrow Z$, and $\cF$ a
quasi--coherent sheaf on $Z$. A {\it $\Gamma$--linearization} of $\cF$
is an action of $\Gamma$ on the symmetric algebra $S(\cF)$ by graded
automorphisms, compatibly with its action on $\cO_Z$. Equivalently,
$\cF$ is endowed with an isomorphism
$$
\varphi:\sigma^*\cF\rightarrow p_2^*\cF
$$
of sheaves on $\Gamma\times Z$, where $p_2:\Gamma\times Z\to Z$
denotes the projection; and, in addition, the following cocycle
condition holds: 
$$
(\mu\times 1_Z)^*\varphi = p_{23}^*\varphi \circ 
(1_{\Gamma}\times\sigma)^*\varphi
$$
as sheaves on $\Gamma\times \Gamma\times Z$, where 
$\mu:\Gamma\times\Gamma\to\Gamma$ denotes the multiplication, and 
$p_{23}:\Gamma\times\Gamma\times Z\to \Gamma\times Z$ the
projection.

The $\Gamma$--linearized sheaves on $Z$ are the objects of a category,
with morphisms being those morphisms of sheaves that are compatible
with the linearizations. Many operations on sheaves extend to
$\Gamma$--linearized sheaves; for example, given two such sheaves
$\cF$, $\cG$ on $Z$, the sheaves $\Ext^n_Z(\cF,\cG)$ are linearized 
as well. Moreover, for any $\Gamma$--module $V$, the sheaf 
$\cF\otimes_k V$ is $\Gamma$--linearized; this yields twists 
$\cF(\chi)$ by characters $\chi:\Gamma\to \mG_m$. If, in addition, 
$Z$ is Cohen--Macaulay, then its dualizing sheaf $\omega_Z$ is
$\Gamma$--linearized as well. 

We will use the following observation: Assume that 
$\Gamma$ decomposes as a product $\Gamma_1\times\Gamma_2$ and that
the $\Gamma_1$-action on $Z$ admits a $\Gamma_2$--equivariant quotient
$q:Z\to Z/\Gamma_1$ which is a locally trivial principal
$\Gamma_1$--bundle. Then, for any $\Gamma$-linearized sheaf $\cF$ on
$Z$, the $\Gamma_1$--invariant direct image $(q_*\cF)^{\Gamma_1}$ is
a $\Gamma_2$--linearized sheaf on $Z/\Gamma_1$, and the canonical map
$q^*(q_*\cF)^{\Gamma_1}\to\cF$ is an isomorphism. On the other hand,
for any $\Gamma_2$--equivariant sheaf $\cG$ on $Z/\Gamma_1$, the sheaf
$q^*\cG$ is $\Gamma$--linearized, and the canonical map 
$\cG\to q_*^{\Gamma_1}(q^*\cG)$ is an isomorphism. Thus, the
categories of $\Gamma$-equivariant sheaves on $Z$ and of 
$\Gamma_1$--equivariant sheaves on $Z/\Gamma_2$ are equivalent.

We now return to a $B$--variety $Y$. Let $P$ be a
standard parabolic subgroup. The group $P\times B$ acts on $P\times Y$
via $(p,b)\cdot (q,y)=(pqb^{-1},by)$; the preceding observation
applies to the quotient map 
$$
q:P\times Y\to P\times ^B Y
$$ 
by the $B$--action, and also to the projection 
$$
p_2:P\times Y\to Y,
$$ 
the quotient by the $P$--action. This yields an equivalence between
the categories of $B$--linearized sheaves on $Y$, and of
$P$--linearized sheaves on $P\times^B Y$. Specifically, we obtain the
following 

\begin{lemma}\label{induction}

\begin{enumerate}

\item
For any $B$--linearized sheaf $\cF$ on $Y$, the sheaf
$$
P\times^B\cF=q_*^B(p_2^*\cF)
$$
on $P\times^B Y$ is $P$-linearized, and $i^*(P\times^B \cF)\cong\cF$ 
as $B$--linearized sheaves (where $i:Y\to P\times^B Y$ denotes the
inclusion).

\item
Conversely, there is a natural isomorphism 
$\cG\cong P\times^B i^*\cG$ for any $P$--linearized sheaf $\cG$ on
$P\times^B Y$. This yields an equivalence between the categories of
$B$--linearized sheaves on $Y$, and of $P$--linearized sheaves on
$P\times^B Y$.

\item
For any $B$--linearized sheaves $\cF$, $\cG$ on $Y$ and for any
integer $n$, we have 
$$
\Ext^n_{P\times^B Y}(P\times^B\cF,P\times^B\cG)\cong
P\times^B \Ext^n_Y(\cF,\cG).
$$

\item
$\cO_{P\times^B Y}=P\times^B\cO_Y$ as $P$--linearized sheaves.
If, in addition, $Y$ is Cohen--Macaulay, then so is $P\times^B Y$, and
$$
\omega_{P\times^B Y}=P\times^B \omega_Y(\chi)
$$
as $P$--linearized sheaves, where $\chi$ is the character of the
$B$--action on the fiber of $\omega_{P/B}$ at the base point $B/B$. In
particular, if $P$ is a minimal parabolic subgroup with associated 
simple root $\alpha$, then
$$
\omega_{P\times^B Y}=P\times^B \omega_Y(\alpha).
$$
\end{enumerate}

\end{lemma}

\begin{proof}
(1) Since $p_2$ is $B$--equivariant and $P$--invariant, the sheaf
$p_2^*\cF$ is $P\times B$--linearized. Thus, $P\times^B \cF$ is
$P$--linearized. 

By construction, the $B$--linearized sheaf
$i^*(P\times^B\cF)$ is the $B$--invariant direct image of the sheaf
$p_2^*\cF$ on $B\times Y$, under the quotient map 
$B\times Y\to B\times^B Y\cong Y$. The latter identifies to the
action map $\sigma:B\times Y\to Y$, so that 
$$
i^*(P\times^B\cF) \cong \sigma_*^B(p_2^*\cF)
\cong \sigma_*^B(\sigma^*\cF) \cong \cF.
$$

(2) Note that $q^*\cG$ is a $P\times B$--linearized sheaf on 
$P\times Y$. Since $p_2$ is the quotient by the $P$--action, we have
$q^*\cG = p_2^*\cF$ for a unique $B$--linearized sheaf $\cF$ on $Y$. It
follows that $\cG \cong q_*^B(q^*\cG) = P\times^B\cF$, so that 
$i^*\cG \cong \cF$ by (1). 

(3) Using the local triviality of $p_2$ and $q$, we obtain isomorphisms 
$$\displaylines{
q^*\Ext^n_{P\times^B Y}(P\times^B \cF,P\times^B \cG) \cong 
\Ext^n_{P\times Y}(q^*(P\times^B \cF),q^*(P\times^B \cG))
\hfill\cr\hfill
\cong \Ext^n_{P\times Y}(p_2^*\cF,p_2^*\cG) \cong
p_2^*\Ext^n_Y(\cF,\cG).
\cr}$$
Now applying $q_*^B$ yields our result.

(4) We have $P\times^B\cO_Y=q_*^B\cO_{P\times Y} \cong \cO_{P\times^B Y}$. 

Recall that the map $f:P\times^B Y\to P/B$ is locally trivial, with
fiber $Y$ at the base point $B/B$. Thus, if $Y$ is Cohen--Macaulay,
then so is $P\times^B Y$. Then the sheaf $\omega_{P\times^B Y}$ is
$P$--linearized, and we obtain
$i^*\omega_{P\times^B Y} \cong \omega_Y(\chi)$. By (1) and (2), it 
follows that $\omega_{P\times^B Y} \cong P\times^B \omega_Y(\chi)$.
\end{proof}

We will also need the following easy result, proved e.g. in
\cite{Bn01} p.~294.

\begin{lemma}\label{vanishing}
Let $Y$ be a $B$-variety and $P$ a minimal parabolic subgroup. Then:

\begin{enumerate}

\item
$R^n\pi_*\cF=0$ for any $n\ge 2$ and any coherent sheaf $\cF$ on
$P\times^B Y$.

\item
$R^1\pi_*\cO_{P\times^B Y}=0$.
\end{enumerate}

\end{lemma}

%

As a final preparation for the proof of Theorem \ref{main}, we
formulate a slightly stronger version of \cite{Bn01} Lemma 8, which
is proved by the same argument.

\begin{lemma}\label{cancellation}
Let $\cF$ a coherent $B$--linearized sheaf on $Y$ satisfying the
following assumptions.

\begin{enumerate}

\item
$\Supp\cF$ is $N_G(Y)$--invariant.

\item
$\pi_*(P\times^B\cF)=0$ for any minimal parabolic subgroup $P$ raising
$Y$.

\end{enumerate}

Then $\Supp\cF$ is $G$--invariant.
\end{lemma}

%

\section{Proof of Theorem \ref{main}}

We may assume that $G\cdot Y=X$ is normal. We show that $Y$ is normal
by induction on its codimension in $X$. Thus, we assume that $Y\ne X$,
and that $P\cdot Y$ is normal for any minimal parabolic subgroup $P$
raising $Y$.

Let $\nu:Z\to Y$ be the normalization, then we have a short exact
sequence 
$$
0\rightarrow \cO_Y \rightarrow \nu_*\cO_Z \rightarrow \cF 
\rightarrow 0\leqno (1)
$$
of sheaves on $Y$, where $\Supp\cF$ is the non--normal
locus. We will show that $\cF=0$.

The normalizer $N_G(Y)$ acts on $Z$, and $\nu$ is equivariant. 
Therefore, (1) is an exact sequence of $N_G(Y)$--linearized sheaves.
Let $P$ be a minimal parabolic subgroup, then we obtain an exact sequence
$$
0 \rightarrow P\times^B\cO_Y 
\rightarrow P\times^B\nu_*\cO_Z 
\rightarrow P\times^B \cF 
\rightarrow 0
$$
of $P$--linearized coherent sheaves on $P\times^B Y$, whence an exact
sequence 
$$
0 \rightarrow \pi_*(P\times^B\cO_Y) 
\rightarrow \pi_*(P\times^B\nu_*\cO_Z) 
\rightarrow \pi_*(P\times^B \cF) 
\rightarrow R^1\pi_*(P\times^B\cO_Y)
$$
of $P$--linearized sheaves on $P\cdot Y$. Now the map 
$$
P\times^B\nu : P\times^B Z \rightarrow P\times^B Y
$$
is the normalization, and 
$P\times^B(\nu_*\cO_Z)=(P\times^B \nu)_*\cO_{P\times^B Z}$. 
Moreover, $P\times^B\cO_Y\cong \cO_{P\times^B Y}$ by Lemma
\ref{induction}, and $R^1\pi_*\cO_{P\times^B Y}=0$ by Lemma
\ref{vanishing}. This yields an exact sequence
$$
0 \rightarrow \pi_*\cO_{P\times^B Y} \rightarrow
\pi_*(P\times^B \nu)_*\cO_{P\times^B Z}
\rightarrow \pi_*(P\times^B\cF) \rightarrow 0.\leqno(2)
$$
If $P$ raises $Y$, then $P\cdot Y$ is normal and both morphisms 
$\pi$, $\pi\circ (P\times^B\nu)$ are proper and birational. 
Thus, by Zariski's main theorem, the maps 
$\cO_{P\cdot Y}\to \pi_*\cO_{P\times^B Y}$
and $\cO_{P\cdot Y}\to\pi_*(P\times^B \nu)_*\cO_{P\times^B Z}$ 
are isomorphisms. Thus, (2) implies $\pi_*(P\times^B\cF)=0$. By 
Lemma \ref{cancellation}, we conclude that $\cF=0$.

Next we assume, in addition, that $X$ is Cohen--Macaulay, and we
show that $Y$ is Cohen--Macaulay as well. For this, we argue 
again by induction on $\codim_X(Y)=c$: we assume that 
$P\cdot Y$ is Cohen--Macaulay for any minimal parabolic subgroup $P$
raising $Y$. We will show that $\Ext^n_X(\cO_Y,\omega_X)=0$ for all
$n\ne c$.

Let $P$ be a minimal parabolic subgroup raising $Y$. The canonical
morphism
$$
\varphi:P\times^B X \rightarrow X
$$ 
identifies with the projection $p_2:P/B\times X\to X$, with fibers
$\mP^1$. By duality for this proper morphism, we obtain an
isomorphism
$$
R\varphi_*R\Hom_{P\times^B X}(\cO_{P\times^B Y},\omega_{P\times^B X}[1])
\cong 
R\Hom_X(R\varphi_*\cO_{P\times^B Y},\omega_X).\leqno(3)
$$
Now $\varphi_*\cO_{P\times^B Y}=\cO_{P\cdot Y}$ by normality of $Y$,
and $R^n\varphi_*\cO_{P\times^B Y}=0$ for all $n\geq 1$, by Lemma
\ref{vanishing}. Thus, (3) yields a spectral sequence 
$$
R^p\varphi_*\Ext^q_{P\times^B X}
(\cO_{P\times^B Y},\omega_{P\times^B X})\Rightarrow
\Ext^{p+q-1}_X(\cO_{P\cdot Y},\omega_X).
$$
But $R^p\varphi_*\cF=0$ for any $p\ge 2$ and any coherent sheaf
$\cF$, by Lemma \ref{vanishing}. Thus, the spectral sequence (3)
reduces to short exact sequences
$$
\displaylines{
(4) \quad \quad 0\rightarrow \varphi_*\Ext^n_{P\times^B X}
(\cO_{P\times^B Y},\omega_{P\times^B X}) \rightarrow
\Ext^{n-1}_X(\cO_{P\cdot Y},\omega_X) 
\hfill\cr\hfill
\rightarrow R^1\varphi_*\Ext^{n-1}_{P\times^B X}
(\cO_{P\times^B Y},\omega_{P\times^B X}) \rightarrow 0.
\cr}$$
Moreover, $\Ext^{n-1}_X(\cO_{P\cdot Y},\omega_X)=0$ for all $n\ne c$,
since both $X$ and $P\cdot Y$ are Cohen--Macaulay and 
$\codim_X(P\cdot Y)=c-1$. This implies
$$
\varphi_*\Ext^n_{P\times^B X}
(\cO_{P\times^B Y},\omega_{P\times^B X})=0
$$
for all $n\ne c$. Together with Lemma \ref{induction}, it follows that
$$
\varphi_*(P\times^B \Ext^n_X(\cO_Y,\omega_X(\alpha)))=0
$$
for all $n\ne c$. But
$\Ext^n_X(\cO_Y,\omega_X(\alpha))\cong 
\Ext^n_X(\cO_Y,\omega_X)(\alpha)$
as a $B$--linearized sheaf on $X$, and this sheaf is killed by
the ideal sheaf $\cI_Y$. This yields for $n\ne c$:
$$
\pi_*(P\times^B \Ext^n_X(\cO_Y,\omega_X)(\alpha))=0,
$$
Moreover, the sheaf $\Ext^n_X(\cO_Y,\omega_X)$ is
$N_G(Y)$--linearized, so that 
$$
\Supp\Ext^n_X(\cO_Y,\omega_X)(\alpha)=\Supp\Ext^n_X(\cO_Y,\omega_X)
$$
is $N_G(Y)$--invariant. By Lemma \ref{cancellation}, it follows that 
$\Ext^n_X(\cO_Y,\omega_X)=0$ for $n\ne c$.

\medskip

\noindent
{\bf Remark 1.} Theorem \ref{main} yields a direct proof for
normality and Cohen--Macaulayness of Schubert varieties. It also
applies to the {\it large Schubert varieties} of \cite{BP00}. 

Specifically, let $G$ be adjoint, i.e., with trivial center. Then $G$
admits a canonical $G\times G$--equivariant completion $X$, a
nonsingular projective variety containing only finitely many $B\times
B$--orbits; moreover, the $G\times G$--orbit closures in $X$ are
nonsingular. Let $Y\subseteq X$ be a $B\times B$--orbit closure, then
$Y$ is multiplicity--free and contains no $G\times G$--orbit, as
follows from \cite{Sr02} \S 2. Since $(G\times G)\cdot Y$ is
nonsingular, Theorem 2 yields that $Y$ is normal and
Cohen--Macaulay. 

This was proved in \cite{BP00} for the large Schubert varieties, the
closures in $X$ of the $B\times B$--orbits in $G$, by a more
complicated argument.

\medskip

\noindent
{\bf Remark 2.} Recall that a morphism $\pi:\tZ\to Z$ of varieties is
{\it rational} if it satisfies the following conditions:

\begin{enumerate}

\item 
$\tZ$ is normal and Cohen--Macaulay.

\item
$\pi$ is proper and birational.

\item
$\pi_*\cO_{\tZ}=\cO_Z$ and $R^n\pi_*\cO_{\tZ}=R^n\pi_*\omega_{\tZ}=0$
for all $n\ge 1$.

\end{enumerate}

Then $Z$ is normal and Cohen--Macaulay, and 
$\pi_*\omega_{\tZ}=\omega_Z$, as follows from duality for the morphism
$\pi$. 

Now, with the notation of Theorem \ref{main}, the morphism 
$\pi_Y:P\times^B Y\to P\cdot Y$ is rational, if $G\cdot Y$ is normal
and Cohen--Macaulay. Indeed, normality of $Y$ implies 
$\pi_*\cO_{P\times^B Y}=\cO_{P\cdot Y}$. On the other hand,
$R^n\pi_*\cO_{P\times^B Y}=0$ for all $n\ge 1$, and
$R^n\pi_*\omega_{P\times^B Y}=0$ for all $n\ge 2$, by Lemma 
\ref{vanishing}. In addition, since $Y$ and $X=G\cdot Y$ are
Cohen--Macaulay, we have 
$$
\Ext^c_{P\times^B X}(\cO_{P\times^B Y},\omega_{P\times^B X})
=\omega_{P\times^B X},
$$ 
and $\Ext^n_{P\times^B X}(\cO_{P\times^B Y},\omega_{P\times^B X})=0$
for all $n\ne c$. 
Thus, the exact sequence (4) with $n=c+1$ yields
$R^1\pi_*\omega_{P\times^B Y}=0$ (and this exact sequence with $n=c$
yields an isomorphism
$\pi_*\omega_{P\times^B Y}\cong\omega_{P\cdot Y}$).

Moreover, given an equivariant morphism 
$$
\pi: \tZ \to Z
$$ of 
$B$--varieties, the induced morphism 
$$
P\times^B\pi : P\times^B\tZ \to P\times^B Z
$$ 
is rational as well, for any standard parabolic subgroup $P$. As a
consequence, the morphism 
$$
P_1\times^B P_2\times^B \cdots \times^B P_d\times^B Y
\rightarrow 
P_1 P_2\cdots P_d\cdot Y
$$
is rational, for any sequence $(P_1,P_2,\ldots,P_d)$ of minimal
parabolic subgroups such that each $P_i$ raises 
$P_{i+1}\cdots P_d\cdot Y$. 

When applied to Schubert varieties, this yields an alternative proof for
rationality of their standard resolutions.

\medskip

\noindent
{\bf Remark 3.} Assume that the characteristic of $k$ is zero. Recall
that a variety $Z$ has {\it rational singularities} if it admits a
rational desingularization; then every desingularization of $Z$ is
rational. Now let $X$ be a $G$--variety and $Y\subseteq X$ a
multiplicity--free $B$--subvariety containing no $G$--orbit ; if
$G\cdot Y$ has rational singularities, then the same holds for $Y$.

This follows indeed from the proof of \cite{Bn01} Theorem 5.
The assumption that $Y$ contains no $G$--orbit was omitted in the
statement of that Theorem, but in fact it is used in its proof, as 
the argument on p.~295, l.~22--25 that would allow to get rid of this 
assumption, is incorrect (I thank Prof.~G.~Zwara for
pointing out this error to me). As a consequence, Theorem 5 in
[loc.cit.] is only proved for regular varieties. This does not
affect the subsequent results of [loc.cit.], since they only concern
regular varieties.

\section{Proof of Theorem \ref{cone}}

By Theorem \ref{main}, $V$ is normal and Cohen--Macaulay. We now
construct a flat family of subvarieties of $\Fl$ with general fiber
$V$ and special fiber a reduced Cohen--Macaulay union of Schubert
varieties.

Write $\Fl=G/P$. Regard $V$ as a point of the Chow variety $\Chow(G/P)$ 
parametrizing all effective cycles in $G/P$ of appropriate
dimension and degree. Recall that $\Chow(G/P)$ is a projective 
$G$--variety; thus, by Borel's fixed point theorem, the orbit closure 
$\overline{B\cdot V}$ contains $B$--fixed points. Choose such a fixed
point $V_0$. Then there exist a nonsingular irreducible curve $C$, a
closed point $c_0\in C$, and a morphism $\varphi: C\to
\overline{B\cdot V}$ such that: $\varphi(c_0)=V_0$, and 
$\varphi(C\setminus \{c_0\})\subseteq B\cdot V$. Since the variety $B$
is rational, we may further assume that $C$ is rational as well
(\cite{Fn98} Example 10.1.6). In particular, $V$ is rationally
equivalent to the cycle $V_0$. Since the latter is $B$--invariant,
it is a multiplicity--free sum of Schubert cycles.

Let $\cV$ be the closure in $G/P\times C$ of the subset 
$$
\{(x,c)\in G/P\times (C\setminus\{c_0\}) ~\vert~ x\in\varphi(c)\}.
$$
Let $p:\cV \to C$ be the restriction of the projection 
$G/P\times C\to C$. Then $\cV$ is a variety, and the morphism $p$ 
is projective and flat. Its fibers at all closed points of
$C\setminus\{c_0\}$ are isomorphic to $V$, whereas the support of 
the fiber $p^{-1}(c_0)$ is contained in the support of $V_0$.
But by flatness of $p$, the associated cycle of $p^{-1}(c_0)$ 
is rationally equivalent to $V$, and hence to $V_0$. Now the linear 
independence of Schubert cycles implies the equality of this 
associated cycle with $V_0$. 

Next let $\cY=\{(g,c)\in G\times C~\vert~ (g^{-1}P/P,c)\in \cV\}$. 
This is a subvariety of $G\times C$; the latter is a $G$--variety 
(for the action by left multiplication on $G$), and $\cY$ is 
$B$--invariant. Clearly, $G\cdot \cY =G\times C$, and $\cY$ contains 
no $G$--orbit. Moreover, $\cY$ is multiplicity-free, as follows from
Lemma \ref{multfree}. Thus, $\cY$ is normal and Cohen--Macaulay by
Theorem \ref{main}, so that $\cV$ is normal and Cohen--Macaulay as
well, by Lemma \ref{multfree} again. It follows that the fiber
$p^{-1}(c_0)$ is equidimensional and Cohen--Macaulay. On the other
hand, this fiber is generically reduced, by the preceding step. Thus,
$p^{-1}(c_0)$ is reduced, and we may identify it to $V_0$.

So $p$ realizes a flat degeneration of $V$ to $V_0$, a reduced
Cohen--Macaulay union of Schubert varieties. Now the
restriction map $H^0(G/P,\cL)\to H^0(V_0,\cL)$ is surjective for any
globally generated invertible sheaf $\cL$ on $G/P$, and
$H^n(V_0,\cL)=0$ for any $n\ge 1$, by \cite{Rn87} 
(1.9, 1.13, 3.5, 3.7). By semicontinuity, the same holds for $V$. 

If, in addition, $\cL$ is ample, then $H^n(V_0,\cL^{-\nu})=0$ for all
$n\le \dim(V_0)-1=\dim(V)-1$ and $\nu\gg 0$, since $V_0$ is
equidimensional and Cohen--Macaulay. But $V_0$ is Frobenius split in
positive characteristics, so that $H^n(V_0,\cL^{-1})=0$ for all
$n\le\dim(V)-1$ by \cite{Rn87} 1.12, 3.7. By semicontinuity again, the
same holds for $V$. 

By \cite{Rn85} Theorem 5, the flag variety $G/P$ is arithmetically
normal in the projective embedding given by any ample invertible sheaf
$\cL$. Since $V$ is normal and the restriction map 
$H^0(G/P,\cL)\to H^0(V,\cL)$ is surjective, it follows that $V$ is
arithmetically normal as well. In addition, $V$ is Cohen--Macaulay and 
$H^n(V,\cL^{\nu})=0$ for $1\le n\le \dim(V)-1$ and all integers $\nu$, 
so that $V$ is arithmetically Cohen--Macaulay.

\end{document}